\documentclass[1p]{elsarticle}
\usepackage{amssymb}
\usepackage{amsfonts}

\usepackage[polish, english]{babel}
\usepackage{polski}

\newtheorem{theorem}{Theorem}
\newtheorem{conjecture}[theorem]{Conjecture}
\newtheorem{corollary}[theorem]{Corollary}
\newtheorem{lemma}[theorem]{Lemma}
\newtheorem{observation}[theorem]{Observation}
\newtheorem{problem}[theorem]{Problem}

\newproof{pf}{Proof}

\begin{document}

\title{A note on the conflict-free chromatic index}

\author[agh]{Mateusz Kamyczura}
\author[agh]{Mariusz Meszka}
\author[agh]{Jakub Przyby{\l}o}

\address[agh]{AGH University of Science and Technology, Faculty of Applied Mathematics, al. A. Mickiewicza 30, 30-059 Krakow, Poland}

\begin{abstract}
Let $G$ be a graph with maximum degree $\Delta$ and without isolated vertices. An edge colouring $c$ of $G$ is conflict-free if the closed neighbourhood of every edge includes a uniquely coloured element. The least number of colours admitting such $c$ is the conflict-free chromatic index of $G$, denoted by $\chi'_{CF}(G)$. In ``\emph{Conflict-free chromatic number versus conflict-free chromatic index}'' [J. Graph Theory, 2022; 99: 349--358] it was recently proved by means of the probabilistic method that $\chi'_{CF}(G)\leq C_1\log_2\Delta+C_2$, where $C_1>337$ and $C_2$ are constants, whereas there are families of graphs with $\chi'_{CF}(G)\geq (1-o(1))\log_2\Delta$.  In this note we provide an explicit simple proof of the fact that $\chi'_{CF}(G)\leq 3\log_2\Delta+1$, which is a corollary of a stronger result: $\chi'_{CF}(G)\leq 3\log_2\chi(G)+1$. For this aim we prove a few auxiliary observations, implying in particular that $\chi'_{CF}(G)\leq 4$ for bipartite graphs.
\end{abstract}


\maketitle

\section{Introduction}
Let $G=(V,E)$ be a graph with maximum degree $\Delta$ and minimum degree $\delta>0$. Consider an edge $C$-colouring $c$ of $G$, i.e. an assignment of elements from a set $C$ of potential colours to the edges of $G$, which does not need to be proper. Let $E_G(v)$ denote the set of edges incident to a vertex $v$ in $G$ and set $E_G[uv]:=E_G(u)\cup E_G(v)$ for every edge $uv\in E$. Let further  $P_c(uv)$, $P_c(v)$ be the \emph{palettes} of colours of $uv$ and $v$, i.e. the multisets of colours in $E(uv)$ and $E[v]$, respectively. We say that $uv$, resp. $v$, is \emph{satisfied} by $c$ if there is a unique colour in its palette, 
i.e. a colour associated to exactly one edge in $E[uv]$, rsp. $E(v)$. For convenience we shall be writing up each multiset in the following form: $c_1^{m_1}\ldots c_p^{m_p}$ where $c_i\in C$ and $m_i$ indicates the number of repetitions of $c_i$ in the multiset, and we shall be abbreviating  $c_i^1$ as $c_i$. An edge $e$ is thus in particular satisfied by $c$ if $P_c(e)$ is of the form $c_1c_2^{m_2}\ldots c_p^{m_p}$. The colouring $c$ is said to be \emph{conflict-free} if all edges of $G$ are satisfied, i.e. contain unique colours in their closed neighbourhoods. The least number of colours admitting such a colouring is called the \emph{conflict-free chromatic index} of $G$ and denoted $\chi'_{CF}(G)$. 
By definition, $\chi'_{CF}(G)\leq \chi'(G)$, but $\chi'_{CF}(G)$ is usually much smaller than $\chi'(G)$.
Recently D\k{e}bski and Przyby{\l}o~\cite{DebskiPrzybylo} exhibited that $\chi'_{CF}\geq (1-o(1))\log_2\Delta$ for the family of complete graphs and used the probabilistic method to prove that there is a constant $C_2$ such that 
$\chi'_{CF}(G)\leq 9\log_{\frac{1}{1-e^{-4}}}\Delta+C_2 = (9/\log_2\frac{1}{1-e^{-4}})\log_2\Delta+C_2$,
where $9/\log_2\frac{1}{1-e^{-4}}\approx 337.5$. The main goal of the present note is to provide a simple direct argument implying that $\chi'_{CF}(G)\leq 3\log_2\Delta+1$. In fact we shall prove something even stronger, 
relating the graph invariant under investigation with the chromatic number $\chi(G)$ rather than $\Delta(G)$,
namely that $\chi'_{CF}(G)\leq 3\log_2\chi(G)+1$, cf. Corollaries~\ref{CorollaryGeneralChi} and~\ref{CorollaryGeneralDelta}. For this purpose we study bipartite graphs first and a few related issues. 

The problem above grew out of its earlier variant concerning vertex colourings. The naturally defined graph parameter in vertex setting is called the \emph{conflict-free chromatic number} and denoted $\chi_{CF}(G)$. In fact investigating $\chi'_{CF}(G)$ goes down to study of $\chi_{CF}(G)$ for the family of line graphs, cf.~\cite{DebskiPrzybylo}. Originally the concept of $\chi_{CF}(G)$ was motivated by channel assignment in wireless networks, where colours represent available frequencies, which potentially interfere, and thus one should always have in range a transmitter with a unique frequency associated, see e.g.~\cite{EvenEtAl,SmorodinskyPhd,SmorodinskyApplications}. 
In 2009 Pach and Tardos~\cite{PachTardos} showed that $\chi_{CF}(G) = O(\ln^{2+\epsilon}\Delta)$. 
This was later improved to $\chi_{CF}(G) = O(\ln^2\Delta)$ by Bhyravarapu, Kalyanasundaram and Mathew~\cite{Hindusi},
while Glebov, Szab\'o and Tardos~\cite{GlebovEtAl} constructed a family of graphs with $\chi_{CF}(G) = \Omega(\ln^2\Delta)$, legitimating the optimal order of the bound from~\cite{Hindusi}. See also~\cite{Hindusi2,KellerEtAl,KostochkaEtAl} for other related results.

\section{Results}

A palette including just one colour shall be called \emph{monochromatic}. 
Moreover, ginen a graph $G=(V,E)$ and an edge $uv$ (which does not have to belong to $E$), we set $G+uv:=(V\cup\{u,v\},E\cup\{uv\})$.

In~\cite{DebskiPrzybylo} it was in particular proved that $\chi'_{CF}(T)\leq 3$ for every tree $T$. 
Here we rephrase this result in order to emphasize some detailed features of the colouring implying this upper bound.
By \emph{levels} of a tree $T$ rooted at $r$ we mean $L_0,L_1,L_2,\ldots$ such that $L_i=\{v\in V(T):{\rm dist}(v,r)=i\}$, where ${\rm dist}(v,r)$ denotes the distance of $v$ and $r$ in $T$. If $d_T(r)=1$, i.e. $r$ is a leaf, we refer to its only incident edge as the \emph{rooted edge}.
A tree of order $1$ is called \emph{trivial}.

\begin{lemma}\label{TreeColouring}
Let $L_0,\ldots,L_q$ be levels of a nontrivial tree $T$ rooted at its leaf $r$.
Then there exists an edge $\mathbb{Z}_3$-colouring $c$ of $T$ such that for every $i\geq 1$ and $u\in L_i$: 
\begin{enumerate} 
\item[($1^\circ$)] $P_c(u)=(i-1)i^{k_u}$, where $k_u\geq 0$ (with $i,i-1$ taken modulo $3$);
\item[($2^\circ$)] if moreover $uv\in E(T)$ where $v\in L_{i+1}$, then $P_c(uv)=(i-1)i^{k_e}(i+1)^{k'_e}$ where $k_e,k'_e\geq 0$. 
\end{enumerate}
Consequently, $\chi'_{CF}(T)\leq 3$. 
\end{lemma}
\begin{pf}
Set $c(e)=i$ (mod $3$) for every edge $e$ between levels $L_i$ and $L_{i+1}$.
It is straightforward to verify ($1^\circ$) for such $c$, while ($2^\circ$) is directly implied by ($1^\circ$).
As $P_c(e)=01^{k_e}$ for the rooted edge $e$, $\chi'_{CF}(T)\leq 3$ follows by ($2^\circ$).
\qed
\end{pf}

We note that it is easy to construct trees for which  $\chi'_{CF}(T)= 3$, consider e.g. the complete binary tree of height $3$.
In order to prove our main result we shall in a way need to extend the upper bound from Lemma~\ref{TreeColouring} towards bipartite graphs in Theorem~\ref{ChiSCF-bipartite}.  
Let us first briefly discuss the base case of complete bipartite graphs, which settle a ``lower bound'' for such a problem.

\begin{observation}\label{ObservComplBip}
If $G=K_{n,m}$ with $n,m\geq 1$, then
$\chi'_{CF}(G)\leq 2$ if $\min\{n,m\}\leq 2$ and $\chi'_{CF}(G)= 3$ if $\min\{n,m\}\geq 3$.
We also have $\chi'_{CF}(C_k)= 2$ for every $k\geq 3$.
\end{observation}
\begin{pf}
Let us assume that $G=K_{n,m}$ with $n,m\geq 3$, as the other case as well as the observation for cycles are straightforward.

A colouring $c'$ implying $\chi'_{CF}(G)\leq 3$ may be constructed by selecting any vertex $v$ in $G$ and painting its incident edges with $0$ and $1$ so that $P_{c'}(v)=01^{d(v)-1}$, and assigning $2$ to all remaining edges.

In order to finally show that $\chi'_{CF}(G)> 2$, suppose on the contrary that $c$ is an edge $\mathbb{Z}_2$-colouring of $K_{n,m}$ satisfying all its edges.
It is easy to see that each colour must be used on at least two edges.
Consider any edge. Since it is satisfied, at least one of its end-vertices, say $v$ must have a unique colour, say $0$ in its palette.
Now consider any other edge coloured $0$ in $G$ and denote by $u$ its end belonging to a different set of the bipartition of $K_{n,m}$ than $v$. Then however $uv$ is incident with at least two edges coloured $0$ and at least two edges coloured $1$, a contradiction.
\qed
\end{pf}

From now on it shall be more convenient for us to focus on a slightly modified graph invariant.
Suppose $c$ is an edge colouring of some subgraph $H$ of a given graph $G=(V,E)$.
Similarly as before we say that an edge $e\in E$ or a vertex $v\in V$ is satisfied by $c$ if $P_c(e)$ or $P_c(v)$, resp., contains a unique colour, where 
$P_c(e)=\{c(e'):e'\in E_H[e]\}$ and $P_c(v)=\{c(e'):e'\in E_H(v)\}$ for $v\in V$.
Then by $\chi'_{sCF}(G)$ we denote the least number of colours admitting an edge colouring $c$ of some subgraph of $G$ (i.e. minimized over all subgraphs) so that all edges of $G$ are satisfied by $c$.
As given such a colouring one may use yet another additional colour to obtain a desired colouring of entire graph, for every graph $G$ without isolated edges we have:
\begin{equation}\label{ChiCFandChiSCF_relation}
\chi'_{sCF}(G)\leq \chi'_{CF}(G)\leq \chi'_{sCF}(G)+1.
\end{equation}
Note that the colouring used in the proof of Observation~\ref{ObservComplBip} actually yields something more.
\begin{observation}\label{ObservComplBipSubgraph} 
For every $n,m\geq 1$, $\chi'_{sCF}(K_{n,m})\leq 2$.
\end{observation}

\begin{theorem}\label{ChiSCF-bipartite}
If $G$ is a bipartite graph without isolated vertices, then $\chi'_{sCF}(G)\leq 3$.
\end{theorem}
\begin{pf}
We may assume $G=(V,E)$ is connected. 
Let $T$ be a BFS spanning tree of $G$ rooted at $r\in V$ with $d_G(r)=\delta(G)$
and let $L_0,L_1,\ldots,L_q$ be the levels of $T$, i.e. $L_0=\{r\}$. Recall that every edge of $G$ must join vertices from two consecutive levels of $T$.

Suppose first that $d_G(r)=1$. Consider a $\mathbb{Z}_3$-colouring $c'$ of $T$ complying with Lemma~\ref{TreeColouring}.
To show that $c'$ satisfies all edges, consider any $uv\in E$ with $u\in L_i,v\in L_{i+1}$ where $i\geq 1$. Then however, analogously as in Lemma~\ref{TreeColouring}, by ($1^\circ$), $P_c(uv)=(i-1)i^{k_e}(i+1)^{k'_e}$ where $k_e,k'_e\geq 0$ (i.e. ($2^\circ$) holds for the edges of $G$). Thus we are done, as the rooted edge is satisfied by $c'$ as well.

Assume therefore that  $d_G(r)=d_T(r)\geq 2$. Let $s$ be the last neighbour of $r$ within the BFS-vertex ordering corresponding to $T$. Let $T_r,T_s$ be the subtrees obtained from $T$ by deleting the edge $rs$ such that $r\in V(T_r)$ and $s\in V(T_s)$. The choice of $s$ shall be crucial in our argument -- note it implies that: 
\begin{eqnarray} 
&&{\rm ~if} ~uv\in E\smallsetminus\{st\} {\rm ~is ~an ~edge ~joining ~levels} ~L_i {\rm ~and} ~L_{i+1} \nonumber\\ 
&&{\rm ~such ~that } ~u\in V(T_r) {\rm ~and} ~v\in (T_s), {\rm ~then} ~u\in L_{i+1} {\rm ~and} ~v\in L_i. \label{CrucialFeature}
\end{eqnarray}
Let $G_r:=G[V(T_r)]$, $G_s:=G[V(T_s)]$, and set $T'_r:=T_r+rs$, $T'_s:=T_s+rs$, $G'_r:=G_r+rs$, $G'_s:=G_s+rs$. Note that all edges of $G$ outside $E(G'_r)\cup E(G'_s)$ must satisfy~(\ref{CrucialFeature}). Note also that $T'_r$ can be regarded as a BFS spanning tree of $G'_r$ rooted at $s$, while $T'_s$ -- as a BFS spanning tree of $G'_s$ rooted at $r$, where $d_{T'_r}(s)=1$, $d_{T'_s}(r)=1$. 
Let $c_r,c_s$ be $\mathbb{Z}_3$-colourings of $T'_r,T'_s$, resp., complying with Lemma~\ref{TreeColouring}, and let $c$ be a concatenation of $c_r$ and $c_s$ (note these agree on $rs$), which is an edge $\mathbb{Z}_3$-colouring of $T$.
Since $T_r$ and $T_s$ are vertex disjoint, then analogously as in the first case, $c$ satisfies all edges in $G_r$ and all edges in $G_s$, while $P_c(rs)=01^{k_{rs}}$ for some $k_{rs}>1$. It thus remains to show that the edges joining $V(T_r)$ and $V(T_s)$, other than $rs$ are satisfied by $c$ as well.

Let $L^r_0,L^r_1,\ldots,L^r_{q_r}$ and $L^s_0,L^s_1,\ldots,L^s_{q_s}$ be the levels in $T'_r$ and $T'_s$, resp., hence in particular $L^r_0=\{s\}$ and $L^s_0=\{r\}$. Note that for every $j\geq 0$:
\begin{eqnarray}\label{LevelsInclusions}
L^s_j\subseteq L_j~~~~~~ {\rm and}~~~~~~L^r_{j+1}\subseteq L_j.
\end{eqnarray}
Suppose $uv\neq rs$ is an edge of $G$ joining levels $L_i$ and $L_{i+1}$, $i\geq 1$, such that $u\in V(T_r)$  and $v\in (T_s)$. By~(\ref{CrucialFeature}), $u\in L_{i+1}$  and $v\in L_i$. Thus by~(\ref{LevelsInclusions}), $u\in L^r_{i+2}$  and $v\in L^s_i$. Therefore, as $c_r,c_s$ comply with Lemma~\ref{TreeColouring}, by ($1^\circ$), $P_c(u)=P_{c_r}(u)=(i+1)(i+2)^{k_u}=(i+1)(i-1)^{k_u}$ (modulo 3) and $P_c(v)=P_{c_s}(v)=(i-1)i^{k_v}$ for some $k_u,k_v\geq 0$. Hence, $P_c(uv)=(i+1)i^{k_{uv}}(i-1)^{k'_{uv}}$ for some $k_{uv},k'_{uv}\geq 0$, and thus $c$ satisfies $uv$.
\qed
\end{pf}

By~(\ref{ChiCFandChiSCF_relation}) we thus obtain the following corollary.

\begin{corollary}\label{ChiCF-bipartite}
If $G$ is a bipartite graph without isolated vertices, then $\chi'_{CF}(G)\leq 4$.
\end{corollary}

We may now prove our main result.

\begin{theorem}\label{ChiSCF-general}
For every graph $G$ without isolated vertices,  $\chi'_{sCF}(G)\leq 3\lceil\log_2\chi(G)\rceil$.
\end{theorem}

\begin{pf}
Let $G=(V,E)$. We prove the theorem by induction. If $\chi(G)=2$, then the assertion follows by Theorem~\ref{ChiSCF-bipartite}.
Let us thus assume that $\chi(G)\geq 3$ and that the assertion holds for graphs with smaller chromatic numbers.
Colour the vertices of $G$ properly with $\chi(G)$ colours and let $V_i$ denote the vertices coloured with $i$th colour, $1\leq i\leq \chi(G)$. Set $V_i=\emptyset$ for $\chi(G)< i\leq 2^{\lceil\log_2\chi(G)\rceil}$, if any; set $k:=2^{\lceil\log_2\chi(G)\rceil-1}$. 
Consider a bipartite subgraph $H$ induced by all edges between $A:=V_1\cup\ldots\cup V_k$ and $B:=V_{k+1}\cup\ldots\cup V_{2k}$ in $G$ 
and let $F$ be the graph induced by all the remaining edges of $G$ (we do not include potential isolated vertices in the vertex sets of $H$ or $F$). Note that $\chi(F)\leq k$, as we may in a way identify colours of the vertices in $V_i$ and $V_{k+i}$ for $1\leq i\leq k$ because there are no edges between these vertices in $F$ (these were included in $F$).

Since $k<\chi(G)$, by induction hypothesis, $\chi'_{sCF}(G)\leq 3\lceil\log_2\chi(F)\rceil \leq 3\lceil\log_2k\rceil = 3\lceil\log_2\chi(G)\rceil - 3$. There thus exists a subgraph of $F$ and its edge colouring with at most $3\lceil\log_2\chi(G)\rceil - 3$ colours which satisfies all edges of $F$. On the other hand, as $H$ is bipartite, by Theorem~\ref{ChiSCF-bipartite} there exists a colouring of a subgraph of $H$ with $3$ additional colours which satisfies all edges in $H$. The concatenation of the both colourings implies the theorem.
\qed
\end{pf}

By~(\ref{ChiCFandChiSCF_relation}), Theorem~\ref{ChiSCF-general} immediately implies the following corollary.
\begin{corollary}\label{CorollaryGeneralChi}
For every graph $G$ without isolated vertices,  $\chi'_{CF}(G)\leq 3\lceil\log_2\chi(G)\rceil+1$.
\end{corollary}

In~\cite{DebskiPrzybylo} it is proved that 
$\chi'_{CF}(K_n) > \lfloor \log_2n-\log_2\log_2n-1\rfloor = \lfloor \log_2\chi(K_n)-\log_2\log_2\chi(K_n)-1\rfloor$.
Thus in general one cannot expect the multiplicative constant of $3$ in the general upper bound from Corollary~\ref{CorollaryGeneralChi} to be pushed down below $1$. Note on the other hand that the same reasoning as applied within the proof of Theorem~\ref{ChiSCF-general} yields a slightly stronger upper bound in the case of the complete graphs.

\begin{observation}\label{UpperBoundCompleteGr}
For every $n\geq 2$, $\chi'_{sCF}(K_n)\leq 2\lceil\log_2 n\rceil$.
\end{observation}
\begin{pf}
The result follows by the same proof as the one of Theorem~\ref{ChiSCF-general}, as for $G=K_n$, $H$ is the complete bipartite graph, and thus by Observation~\ref{ObservComplBipSubgraph}, $\chi'_{sCF}(H)\leq 2$, while $F$ is a union of two disjoint complete graphs.
\qed
\end{pf}

\begin{corollary}\label{CorollaryGeneralDelta}
For every graph $G$ without isolated vertices,  $\chi'_{CF}(G)\leq 3\lceil\log_2\Delta(G)\rceil+1$.
\end{corollary}
\begin{pf}
We may assume $G$ is connected.
Due to Brook's Theorem the result follows directly by Corollary~\ref{CorollaryGeneralChi} for all graphs except cycles and complete graphs. The remaining cases follow by Observations~\ref{ObservComplBip} and~\ref{UpperBoundCompleteGr}.
\qed
\end{pf}

\section{Concluding remarks}

By Corollary~\ref{ChiCF-bipartite} and Observation~\ref{ObservComplBip}  we know that the optimal general upper bound for the conflict-free chromatic index within the family of bipartite graphs is either $3$ or $4$. We propose the following conjecture.
\begin{conjecture}
If $G$ is a bipartite graph without isolated vertices, then $\chi'_{CF}(G)\leq 3$.
\end{conjecture}
Note that even settling this conjecture in the affirmative might not provide 
a significant improvement of the general bound in Corollary~\ref{CorollaryGeneralChi}. This would decrease merely by $1$ on the basis of our approach utilized in the proof of Theorem~\ref{ChiSCF-general}, unless we knew something specific on the corresponding colourings of bipartite graphs, e.g. implying that.$\chi'_{sCF}(G)\leq 2$ for these.

We believe that it would also be interesting to solve the following problem
\begin{problem}
Characterize the family of all trees $T$ with $\chi'_{CF}(T)=3$.
\end{problem}

\end{document}